\documentclass{amsart}

\newtheorem{theorem}{Theorem}[section]
\newtheorem{lemma}[theorem]{Lemma}
\newtheorem{proposition}[theorem]{Proposition}
\theoremstyle{definition}

\newtheorem{example}[theorem]{Example}

\theoremstyle{remark}
\newtheorem{remark}[theorem]{Remark}
\numberwithin{equation}{section}
\newfont{\aj}{eufm10 at10pt}
\newfont{\ajk}{eufm10 at8pt}
\newfont{\kh}{msbm10}
\newcommand{\C}{\mbox{\kh C}}
\newcommand{\cala}{\mbox{\aj A}}
\newcommand{\calb}{\mbox{\aj B}}
\newcommand{\calc}{\mbox{\aj C}}
\newcommand{\calh}{\mbox{\aj H}}
\newcommand{\calx}{\mbox{\aj X}}
\newcommand{\caly}{\mbox{\aj Y}}
\newcommand{\calz}{\mbox{\aj Z}}
\newcommand{\call}{\mbox{\aj L}}
\newcommand{\calk}{\mbox{\aj K}}
\newcommand{\cals}{\mbox{\aj S}}
\newcommand{\calak}{\mbox{\ajk A}}

\newcommand{\calhk}{\mbox{\ajk H}}

\begin{document}

\title{Automatic Continuity of $\sigma$-Derivations on $C^*$-Algebras}
\author{Madjid Mirzavaziri}
\address{Department of Mathematics, Ferdowsi University, P. O. Box 1159, Mashhad 91775, Iran}
\email{mirzavaziri@math.um.ac.ir}
\author{Mohammad Sal Moslehian}
\address{Department of Mathematics, Ferdowsi University, P. O. Box 1159, Mashhad 91775, Iran}
\email{moslehian@ferdowsi.um.ac.ir} \subjclass[2000]{Primary
46L57; Secondary 46L05, 47B47}
\keywords{$*$-$(\sigma,\tau)$-Derivation, $\sigma$-derivation,
derivation, automatic continuity, $C^*$-algebra}

\begin{abstract}
Let ${\calak}$ be a $C^*$-algebra acting on a Hilbert space
${\calhk}$, $\sigma:{\calak}\to B({\calhk})$ be a linear mapping
and $d:{\calak}\to B({\calhk})$ be a $\sigma$-derivation.
Generalizing the celebrated theorem of Sakai, we prove that if
$\sigma$ is a continuous $*$-mapping then $d$ is automatically
continuous. In addition, we show the converse is true in the
sense that if $d$ is a continuous $*$-$\sigma$-derivation then
there exists a continuous linear mapping $\Sigma:{\calak}\to
B({\calhk})$ such that $d$ is $*$-$\Sigma$-derivation. The
continuity of the so-called $*$-$(\sigma,\tau)$-derivations is
also discussed.
\end{abstract}

\maketitle

\section{Introduction.}

Let ${\cala}$ be a subalgebra of an algebra ${\calb}$, ${\calx}$
be a ${\calb}$-bimodule and $\sigma:{\cala}\to {\calb}$ be a
linear mapping. A linear mapping $d:{\cala}\to {\calx}$ is called
a $\sigma$-derivation if $d(ab)=d(a)\sigma(b)+\sigma(a)d(b)$ for
all $a,b\in {\cala}$ (see \cite{M-M} and \cite{MOS} and references
therein). If ${\calb}={\calx}={\cala}$ and $\sigma$ is the
identity map on ${\cala}$, then we reach to the usual notion of a
derivation on the algebra ${\cala}$.

In 1958, I. Kaplansky \cite{KAP} conjectured that every
derivation on a $C^*$-algebra is continuous. H. Sakai \cite{SAK}
proved this conjecture and then J. R. Ringrose \cite{RIN} showed
that every derivation from a $C^*$-algebra ${\cala}$ into a
Banach ${\cala}$-bimodule is continuous.

In this paper, we investigate the continuity of
$\sigma$-derivations on $C^*$-algebras. Let ${\cala}$ be a
$C^*$-algebra acting on a Hilbert space ${\calh}$. We prove that
if $\sigma:{\cala}\to B({\calh})$ is a continuous $*$-linear
mapping then every $\sigma$-derivation from ${\cala}$ into
$B({\calh})$ is automatically continuous and so Sakai's theorem
\cite{SAK} is generalized. In addition, we establish the converse
in the sense that if $d:{\cala}\to B({\calh})$ is a continuous
$*$-$\sigma$-derivation, then there exists a continuous mapping
$\Sigma:{\cala}\to B({\calh})$ such that $d$ is a
$*$-$\Sigma$-derivation. In the last section we discuss the
continuity of the so-called $*$-$(\sigma,\tau)$-derivations.

The importance of our approach is that $\sigma$ is a linear
mapping in general, not necessarily an algebra homomorphism.
There are some applications of $\sigma$-derivations to develop an
approach to deformations of Lie algebras which have many
applications in models of quantum phenomena and in analysis of
complex systems; cf. \cite{H-L-S}.

For the definition and elementary properties of $C^*$-algebras we
refer the reader to \cite{MUR} and \cite{PAL}.

\section{Elementary Properties of $\sigma$-Derivations}

Throughout this section, ${\cala}$ is a subalgebra of an algebra
${\calb}$, ${\calx}$ is a ${\calb}$-bimodule and
$\sigma:{\cala}\to {\calb}$ is a linear mapping.

A linear mapping $d:{\cala}\to {\calx}$ is called a
$\sigma$-derivation if $d(ab)=d(a)\sigma(b)+\sigma(a)d(b)$ for
all $a,b\in {\cala}$. Familiar examples are:

(i) every ordinary derivation $\delta$ of an algebra ${\cala}$
into an ${\cala}$-bimodule ${\calx}$ is an $\iota$-derivation
(here $\iota$ denotes the identity map on ${\cala}$);

(ii) every endomorphism $\varphi$ on ${\cala}$ is a
$\frac{\varphi}{2}$-derivation;

(iii) for a given homomorphism $\sigma$ on ${\cala}$ and a fixed
arbitrary element $x$ in an ${\cala}$-bimodule ${\calx}$, the
so-called $\sigma$-inner derivation is defined to be
$d_x(a)=x\sigma(a)-\sigma(a)x$.

There is an interesting link between $\sigma$-derivations and
algebra homomorphisms as follows.

\begin{theorem} Let $\sigma:{\cala}\to{\calb}$ be a homomorphism and $d:{\cala}\to{\calx}$
be a $\sigma$-derivation. Then

(i) ${\calx}$ equipped with the module multiplications $a\cdot
x=\sigma(a)x$ and $x\cdot a=x\sigma(a)$ is an ${\cala}$-bimodule
denoted by $\widetilde{{\calx}}$;

(ii) $d:{\cala}\to\widetilde{{\calx}}$ is an ordinary derivation;

(iii) ${\calc}={\cala}\oplus\widetilde{{\calx}}$ equipped with the
multiplication $(a,x)(b,y)=(ab,x\cdot b+a\cdot y)$ is an algebra,
and $\varphi_d:{\cala}\to {\calc}$ defined by
$\varphi_d(a)=(a,d(a))$ is an injective homomorphism.

(iv) if ${\cala}, {\calb}$ and ${\calx}$ are normed, $\sigma$ is
continuous, and ${\calc}$ is equipped with the norm
$\|(a,x)\|=\|a\|+\sup\{\|x\|, \|a_1\cdot x\|, \|x\cdot a_2\|,
\|a_1\cdot x\cdot a_2\|: a_1,a_2\in {\cala}, \|a_1\|\leq 1,
\|a_2\|\leq 1\}$, then $\varphi_d$ is continuous if and only if
$d$ is continuous. Thus if every injective homomorphism of
${\cala}$ into a Banach algebra is continuous, then every
$\sigma$-derivation of ${\cala}$ into a Banach ${\calb}$-bimodule
is continuous.
\end{theorem}

\begin{proof} Straightforward (see \cite{PAL}).\end{proof}

Recall that if ${\caly}$ and ${\calz}$ are normed spaces and
$T:{\caly}\to {\calz}$ is a linear mapping, then the set of all
$z$ such that there is a sequence $\{y_n\}$ in ${\caly}$ with
$y_n\to 0$ and $Ty_n\to z$ is called the separating space
${\cals}(T)$ of $T$. Clearly,
${\cals}(T)=\displaystyle{\cap_{n=1}^\infty}\overline{\{T(y):
\|y\|<1/n\}}$ is a closed linear space. If ${\caly}$ and
${\calz}$ are Banach spaces, by the closed graph theorem, $T$ is
continuous if and only if ${\cals}(T)=\{0\}$.

\begin{lemma} Let $d:{\cala}\to {\calx}$ be a $\sigma$-derivation. Then
\begin{eqnarray*}
d(c)(\sigma(ab)-\sigma(a)\sigma(b))=(\sigma(ca)-\sigma(c)\sigma(a))d(b)
\end{eqnarray*}
for all $a, b, c\in{\cala}$.\end{lemma}

\begin{proof}
\begin{eqnarray*}
d(cab)&=&d(c)\sigma(ab)+\sigma(c)d(ab)\\
d(c)\sigma(ab)&=&d(cab)-\sigma(c)d(ab)\\
&=&(d(ca)\sigma(b)+\sigma(ca)d(b))-\sigma(c)d(ab)\\
&=&(d(c)\sigma(a)+\sigma(c)d(a))\sigma(b)+\sigma(ca)d(b)-\sigma(c)d(ab)\\
&=&d(c)\sigma(a)\sigma(b)+\sigma(c)d(a)\sigma(b)+\sigma(ca)d(b)-\sigma(c)d(ab)\\
&=&d(c)\sigma(a)\sigma(b)+\sigma(c)d(a)\sigma(b)+\sigma(ca)d(b)\\
&&-\sigma(c)(d(a)\sigma(b)+\sigma(a)d(b)),
\end{eqnarray*}
whence
\begin{eqnarray*}
d(c)(\sigma(ab)-\sigma(a)\sigma(b))=(\sigma(ca)-\sigma(c)\sigma(a))d(b).
\end{eqnarray*}
\end{proof}
\begin{lemma} Let ${\cala}$ and ${\calb}$ be normed algebras, $\sigma:{\cala}\to{\calb}$ be a continuous mapping and let $d:{\cala}\to{\calb}$ be a $\sigma$-derivation. Then for each $a\in{\cals}(d)$ and $b,c\in{\cala}$ we have $a(\sigma(bc)-\sigma(b)\sigma(c))=0.$
\end{lemma}

\begin{proof} For each $a\in{\cals}(d)$ there exists a sequence
$\{a_n\}$ such that $a_n\to 0$ and $d(a_n)\to a$. By Lemma 2.2 we
have
\[d(a_n)(\sigma(bc)-\sigma(b)\sigma(c))=(\sigma(a_nb)-\sigma(a_n)\sigma(b))d(c)\to
0.\] Thus $a(\sigma(bc)-\sigma(b)\sigma(c))=0$.\end{proof}

\begin{remark} Recall that if $E$ is a subset of an algebra ${\calb}$, the right annihilator $ran(E)$ of $E$ (resp. the left annihilator $lan(E)$ of $E$) is defined to be $\{b\in{\calb}: Eb=\{0\}\}$ (resp. $\{b\in{\calb}: bE=\{0\}\}$). The set $ann(E):=ran(E)\cap lan(E)$ is called the annihilator of $E$. The previous lemma shows that if ${\cala}$ and ${\calb}$ are Banach algebras, $\sigma:{\cala}\to {\calb}$ is a continuous linear mapping, $d:{\cala}\to {\calb}$ is a $\sigma$-derivation and $ran({\cals}(d))=\{0\}$ then $\sigma$ is an endomorphism, and if $lan(\{\sigma(bc)-\sigma(b)\sigma(c):b,c\in{\cala}\})=\{0\}$ then $d$ is continuous.\end{remark}

\begin{proposition} Suppose that ${\cala}$ is a  Banach algebra, ${\calb}$ is a simple Banach algebra, $\sigma:{\cala}\to{\calb}$ is a surjective continuous linear mapping and $d:{\cala}\to{\calb}$ is a $\sigma$-derivation. Then $d$ is continuous or $\sigma$ is an endomorphism.\end{proposition}

\begin{proof} At first we show that ${\cals}(d)$ is a closed bi-ideal
of ${\calb}$. To see this let $b\in{\calb}$ and $a\in{\cals}(d)$.
Then there is a sequence $\{a_n\}$ such that $a_n\to 0$ and
$d(a_n)\to a$. Since $\sigma$ is surjective, there is an element
$c\in{\cala}$ such that $b=\sigma(c)$. Now we have $ca_n\to 0$
and $d(ca_n)=d(c)\sigma(a_n)+\sigma(c)d(a_n)\to 0+ba=ba$. This
shows that $ba\in{\cals}(d)$. By the same way $ab\in{\cals}(d)$.
Thus ${\cals}(d)$ is a bi-ideal.

${\cals}(d)$ is $\{0\}$ or ${\calb}$. If ${\cals}(d)=\{0\}$ then
$d$ is continuous and if ${\cals}(d)={\calb}$ then ${\cals}(d)$
has zero right annihilator and so, by Remark 2.4, $\sigma$ is an
endomorphism.\end{proof}

The following example yields a continuous $\sigma$-derivation
with a non-continuous linear mapping $\sigma$.

\begin{example} Let ${\cala}=C[0,2]$, the $C^*$-algebra of
all complex valued continuous functions defined on the interval
$[0,2]$. Define the continuous function $h:[0,2]\to \C$ by
$h(t)=0$ on $[0,1]$ and $h(t)=t-1$ on $[1,2]$, the linear mapping
$\sigma:{\cala}\to {\cala}$ by
$$\sigma(f)(t)=\left \{ \begin{array}{ll}
\alpha (f|_{[0, 1/2]})(t) & \rm{if~ } 0\le t\le \frac{1}{2} \\
2(1-t) \alpha(f|_{[0,1/2]})(\frac{1}{2})+(t-\frac{1}{2})f(1)
&\rm{if~ } \frac12\le t\le 1\\
\frac{1}{2} f(t)&\rm{if~ } 1\le t\le 2
\end{array}\right.$$
where $\alpha$ is a discontinuous linear mapping on $C[0,1/2]$,
and define the linear mapping $d:{\cala}\to {\cala}$ by
$d(f)=fh$. For all $t\in[0,1]$, we have
\begin{eqnarray*}
d(fg)(t)&=&f(t)g(t)h(t)\\
&=&0\\
&=&f(t)h(t)\sigma(g)(t)+\sigma(f)(t)g(t)h(t)\\
&=&d(f)(t)\sigma(g)(t)+\sigma(f)(t)d(g)(t)\\
&=&(d(f)\sigma(g)+\sigma(f)d(g))(t),
\end{eqnarray*}
and for all $t\in[1,2],$
\begin{eqnarray*}
d(fg)(t)&=&f(t)g(t)h(t)\\&=&f(t)h(t)(\frac{g}2)(t)+(\frac{f}2)(t)g(t)h(t)\\
&=&d(f)(t)\sigma(g)(t)+\sigma(g)(t)d(g)(t)\\
&=&(d(f)\sigma(g)+\sigma(f)d(g))(t).
\end{eqnarray*}
\end{example}

In this example if we define $\Sigma$ on ${\cala}$ by
$\Sigma(f)=\frac{f}2$ then $d$ is a $\Sigma$-derivation and
$\Sigma$ is continuous. In the next section we will show that
this is true in general, i.e, if $d$ is a continuous
$*$-$\sigma$-derivation, then we can find a continuous linear
mapping $\Sigma$ such that $d$ is a $*$-$\Sigma$-derivation.

\section{$\sigma$-Derivations on $C^*$-Algebras}

In this section we establish several significant theorems
concerning the continuity of $\sigma$-derivations on
$C^*$-algebras. Throughout this section, ${\cala}$ denotes a
$C^*$-algebra acting on a Hilbert space ${\calh}$, i.e. a closed
$*$-subalgebra of the algebra $B({\calh})$ consisting of all
bounded linear mappings on ${\calh}$. In addition, we assume that
$\sigma$ and $d$ are linear mappings of ${\cala}$ into
$B({\calh})$.

Our first result states that when we deal with a continuous
$\sigma$-derivation we may assume that $\sigma$ is continuous.

\begin{lemma} Let $d:{\cala}\to{\calb}$ be a continuous $\sigma$-derivation. Then ${\cals}(\sigma)\subseteq ann(d({\cala}))$.\end{lemma}

\begin{proof} Assume $A\in{\cals}(\sigma)$. Then there is a sequence
$\{A_n\}$ in ${\cala}$ such that $A_n\to 0$ and $\sigma(A_n)\to
A$. For each $B\in{\cala}$ we have
$$d(A_nB)=d(A_n)\sigma(B)+\sigma(A_n)d(B)\to d(0)\sigma(B)+Ad(B).$$
Since $d(A_nB)\to d(0)$ and $d(0)=0$, we obtain $Ad(B)=0$.
Similarly, we can show that $d(B)A=0$. \end{proof}

\begin{theorem} Let $\sigma$ be a linear mapping and $d$ be a continuous $*$-$\sigma$-derivation. Then there is a continuous linear mapping $\Sigma:{\cala}\to B({\calh})$ such that $d$ is a $*$-$\Sigma$-derivation.\end{theorem}

\begin{proof} Let ${\call}_0=\bigcup_{A\in{\calak}}d(A)({\calh})$ and
${\call}$ be the closed linear span of ${\call}_0$. Then we can
write ${\calh}={\call}\oplus{\calk}$, where
${\calk}={\call}^\perp$. Thus $0=\langle d(A^*)(h),k\rangle
=\langle h,d(A)k\rangle $ for all $h\in{\calh},k\in{\calk},
A\in{\cala}$ shows that ${\calk}=\bigcap_{A\in{\calak}}\ker d(A)$.
Now define $\Sigma$ on $\cala$ by $\Sigma(A)=\sigma(A)P$ where
$P$ denotes the corresponding projection to ${\call}$.

We show that $d$ is a $\Sigma$-derivation. Let $A,B\in{\cala}$
and $h=\ell+k\in{\call}\oplus{\calk}={\calh}$. Then
\begin{eqnarray*}
d(AB)(h)&=&d(AB)(\ell+k)\\&=&d(AB)(\ell)+d(AB)(k)\\&=&d(AB)(\ell)+0\\
&=&d(A)\sigma(B)(\ell)+\sigma(A)d(B)(\ell)\\&=&d(A)\Sigma(B)(\ell)+\Sigma(A)d(B)(\ell)\\
&=&d(A)\Sigma(B)(\ell)+\Sigma(A)d(B)(\ell)+d(A)\Sigma(B)(k)+\Sigma(A)d(B)(k)\\
&=&d(A)\Sigma(B)(\ell+k)+\Sigma(A)d(B)(\ell+k)\\&=&(d(A)\Sigma(B)+\Sigma(A)d(B))(h).
\end{eqnarray*}

Moreover, $\Sigma$ is continuous on ${\cala}$. Let
$A_n\in{\cala}$ with $A_n\to 0$ and $\Sigma(A_n)\to A$. Then for
each $\ell\in{\call}_0$ there is a $B\in{\cala}$ and there is an
$h\in{\calh}$ such that $\ell=d(B)(h)$. Thus we can write
$A(\ell)=A(d(B)(h))=(Ad(B))(h)=0$, since $A\in ann(d({\cala}))$.
We therefore have $A=0$ on ${\call}_0$ and so $A=0$ on $\call$,
since $A$ is continuous. On the other hand, for each $k\in{\calk}$
we have $0=\Sigma(A_n)(k)\to A(k)$ and so $A=0$ on ${\calk}$. Thus
we have proved that $A=0$ on ${\calh}$. This shows that
${\cals}(\Sigma)=\{0\}$ and so $\Sigma$ is continuous on
${\cala}$.
\end{proof}

\begin{theorem} Let $\sigma$ be a $*$-linear mapping and let $d$ be a continuous $\sigma$-derivation. Then there is a continuous linear mapping $\Sigma:{\cala}\to B({\calh})$ such that $d$ is a $\Sigma$-derivation.\end{theorem}

\begin{proof} Define the $\sigma$-derivation $d^*:{\cala}\to
B({\calh})$ by $d^*(A)=d(A^*)^*$. Let
${\call}_0=(\bigcup_{A\in{\calak}}d(A)({\calh}))\bigcup(\bigcup_{A\in{\calak}}d^*(A)({\calh}))$
and ${\call}$ be the closed linear span of ${\call}_0$. Then we
can write ${\calh}={\call}\oplus{\calk}$, where
${\calk}={\call}^\perp$. Thus $0=\langle d(A^*)(h),k\rangle
=\langle h,d^*(A)k\rangle$ and $0=\langle d^*(A^*)(h),k\rangle
=\langle h,d(A)k\rangle$ for all $h\in{\calh},k\in{\calk},
A\in{\cala}$ shows that ${\calk}=(\bigcap_{A\in{\calak}}\ker
d(A))\bigcap(\bigcap_{A\in{\calak}}\ker d^*(A))$. Now define
$\Sigma$ on $\cala$ by $\Sigma(A)=\sigma(A)P$ where $P$ denotes
the corresponding projection to ${\call}$.

Using the same argument in the proof of Theorem 3.2 one can show
that both $d$ and $d^*$ are $\Sigma$-derivations (note that $d^*$
is a $\sigma$-derivation).

$\Sigma$ is continuous on ${\cala}$. To see this, assume that
$A_n\in{\cala}$, $A_n\to 0$ and $\Sigma(A_n)\to A$. Then for each
$\ell\in{\call}_0$ there is a $B\in{\cala}$ and there is an
$h\in{\calh}$ such that $\ell=d(B)(h)$ or $\ell=d^*(B)(h)$. Thus
we can write $A(\ell)=A(d(B)(h))=(Ad(B))(h)=0$ or
$A(\ell)=A(d^*(B)(h))=(Ad^*(B))(h)=0$, since $A\in
ann(d({\cala}))\cap ann(d^*({\cala}))$. We therefore have $A=0$ on
${\call}_0$ and so $A=0$ on $\call$, since $A$ is continuous. On
the other hand, for each $k\in{\calk}$ we have
$0=\Sigma(A_n)(k)\to A(k)$ and so $A=0$ on ${\calk}$. Thus we
have proved that $A=0$ on ${\calh}$. This shows that
${\cals}(\Sigma)=\{0\}$ and so $\Sigma$ is continuous on
${\cala}$. \end{proof}

The next two propositions allow us to assume that $\sigma$ is a
homomorphism when we discuss the continuity of
$\sigma$-derivations.

\begin{proposition} Let $\sigma$ be a continuous $*$-linear mapping and $d$ be a $\sigma$-derivation. Then there is a continuous $*$-homomorphism $\Sigma:{\cala}\to B({\calh})$ and a $\Sigma$-derivation $D:{\cala}\to B({\calh})$ such that $D$ is continuous if and only if so is $d$. Moreover, if $d$ preserves $*$, then so does $D$.\end{proposition}

\begin{proof} By Lemma 2.3 for each $A\in{\cals}(d)$ and
$B,C\in{\cala}$ we have
\[A(\sigma(BC)-\sigma(B)\sigma(C))=0\hspace{1cm}(*)\]
Now let
${\call}_0=\bigcup_{B,C\in{\calak}}(\sigma(BC)-\sigma(B)\sigma(C))({\calh})$
and $\call$ be the closed linear span of ${\call}_0$. Then $(*)$
implies that $A({\call})=0$ for each $A\in{\cals}(d)$.

We can write ${\calh}={\call}\oplus{\calk}$, where
${\calk}={\call}^\perp$. For each $B,C\in{\cala}, h\in{\calh}$ and
$k\in{\calk}$ we have
\begin{eqnarray*}
0&=&\langle (\sigma(BC)-\sigma(B)\sigma(C))(h),k\rangle \\&=&\langle h,(\sigma(BC)-\sigma(B)\sigma(C))^*(k)\rangle \\
&=&\langle h,(\sigma(C^*B^*)-\sigma(C^*)\sigma(B^*))(k)\rangle .
\end{eqnarray*}
Since $\cala$ is a $*$-subalgebra of $B({\calh})$, we infer that
$(\sigma(BC)-\sigma(B)\sigma(C))(k)=0$ for each $B,C\in{\cala}$
and $k\in{\calk}$. This shows that
${\calk}=\bigcap_{B,C\in{\calak}}\ker(\sigma(BC)-\sigma(B)\sigma(C))$.

Now let $P=P_{\calk}$ be the projection corresponding to
${\calk}$. At first we show that $\sigma(A)P=P\sigma(A)$ for all
$A\in{\cala}$. For each $A,B,C\in{\cala}$ and $k\in{\calk}$ we
have
\begin{eqnarray*}
(\sigma(BC)-\sigma(B)\sigma(C))\sigma(A)(k)&=&(\sigma(BC)\sigma(A)-\sigma(B)\sigma(C)\sigma(A))(k)\\
&=&(\sigma(BCA)-\sigma(BCA))(k)\\&=&0.
\end{eqnarray*}
This shows that $\sigma(A)({\calk})\subseteq{\calk}$ and so
$\sigma(A)P=P\sigma(A)$.

By using Lemma 2.2 we get
\[0=d(B)(\sigma(CA)-\sigma(C)\sigma(A))(k)=(\sigma(BC)-\sigma(B)\sigma(C))d(A)(k)\]
for all $k\in{\calk}$. This implies that $d(A)(k)\in{\calk}$ for
all $k\in{\calk}$. Hence $d(A)({\calk})\subseteq{\calk}$ and so
$d(A)P=Pd(A)$.

Now put $\Sigma(A)=\sigma(A)P$ and $D(A)=d(A)P$ for all
$A\in{\cala}$. Firstly, $\Sigma$ is a $*$-homomorphism. For
$k\in{\calk}$ we have
\begin{eqnarray*}
\Sigma(AB)(k)&=&\sigma(AB)P(k)\\&=&\sigma(AB)(k)\\&=&\sigma(A)\sigma(B)(k)\\
&=&\sigma(A)\sigma(B)P^2(k)\\&=&\sigma(A)P\sigma(B)P(k)\\&=&\Sigma(A)\Sigma(B)(k).
\end{eqnarray*}
And for $\ell\in{\call}$,
\begin{eqnarray*}
\Sigma(AB)(\ell)=\sigma(AB)P(\ell)=0=\sigma(A)P\sigma(B)P(\ell)=\Sigma(A)\Sigma(B)(\ell).
\end{eqnarray*}
Moreover,
$\Sigma(A^*)=\sigma(A^*)P=\sigma(A)^*P=(P\sigma(A))^*=(\sigma(A)P)^*=\Sigma(A)^*
(A\in{\cala})$.

Secondly, $D$ is a $\Sigma$-derivation, since for $A,B\in{\cala}$
and $k\in{\calk}$ we have
\begin{eqnarray*}
D(AB)(k)&=&d(AB)P(k)\\&=&d(AB)(k)\\&=&d(A)\sigma(B)(k)+\sigma(A)d(B)(k)\\
&=&d(A)\sigma(B)P^2(k)+\sigma(A)d(B)P^2(k)\\&=&d(A)P\sigma(B)P(k)+\sigma(A)Pd(B)P(k)\\
&=&(D(A)\Sigma(B)+\Sigma(A)d(B))(k).
\end{eqnarray*}
And for $\ell\in{\call}$,
\begin{eqnarray*}
D(AB)(\ell)&=&d(AB)P(\ell)\\&=&0\\&=&d(A)P\sigma(B)P(\ell)+\sigma(A)Pd(B)P(\ell)\\
&=&(D(A)\Sigma(B)+\Sigma(A)d(B))(\ell).
\end{eqnarray*}

Moreover, $D$ is continuous if and only if so is $d$. To show
this let $D$ be continuous and $A\in{\cals}(d)$. Then there is a
sequence $\{A_n\}$ in ${\cala}$ such that $A_n\to 0$ and
$d(A_n)\to A$. By the first paragraph of the proof we know that
$A({\call})=0$, and for $k\in{\calk}$ we have
$$A(k)=\lim_{n\to\infty}d(A_n)(k)=\lim_{n\to\infty}d(A_n)P(k)=\lim_{n\to\infty}D(A_n)(k)=0,$$
since $D$ is continuous.

On the other hand, if $d$ is continuous and $A\in{\cals}(D)$,
then there is a sequence $\{A_n\}$ in ${\cala}$ such that $A_n\to
0$ and $D(A_n)\to A$. Obviously for $\ell\in{\call}$ we have
\[ A(\ell)=\lim_{n\to\infty}D(A_n)(\ell)=\lim_{n\to\infty}d(A_n)P(\ell)=0,\]
and for $k\in{\calk}$,
$$A(k)=\lim_{n\to\infty}D(A_n)(k)=\lim_{n\to\infty}d(A_n)P(k)=\lim_{n\to\infty}d(A_n)(k)=0,$$
since $d$ is continuous.

Similarly, one can show that $\Sigma$ is continuous. If $d$ is a
$*$-$\sigma$-derivation, then the relation $d(A)P=Pd(A)$ implies
that $D$ is also a $*$-$\Sigma$-derivation. \end{proof}

\begin{remark} If ${\cala}$ is a von Neumann algebra and both $\sigma$ and $d$ are mappings of ${\cala}$ into ${\cala}$ then, due to the fact that the range projection of every element of ${\cala}$ is in ${\cala}$, we conclude that the projection $P$ constructed in the proof of Proposition 3.4 belongs to ${\cala}$ and so the ranges of the derivation $D$ and the linear mapping $\Sigma$ are contained in ${\cala}$.\end{remark}

\begin{proposition} Let $\sigma$ be a continuous linear mapping and $d$ be a $*$-$\sigma$-derivation. Then there is a continuous $*$-homomorphism $\Sigma:{\cala}\to B({\calh})$ and a $*$-$\Sigma$-derivation $D:{\cala}\to B({\calh})$ such that $D$ is continuous if and only if so is $d$.\end{proposition}

\begin{proof} Define $\sigma^*:{\cala}\to B({\calh})$ by
$\sigma^*(A)=\sigma(A^*)^*$. Then $d$ is a $\sigma^*$-derivation.
Clearly $d$ is a $\tau$-derivation where
$\tau=\frac{\sigma+\sigma^*}{2}$ is a continuous $*$-linear
mapping. By Proposition 3.4 there exist a continuous
$*$-homomorphism $\Sigma$ and a $*$-$\Sigma$-derivation $D$ such
that $D$ is continuous if and only if so is $d$. \end{proof}

The following theorem is an extension of Sakai's theorem
\cite{SAK} to $\sigma$-derivations.

\begin{theorem} Let $\sigma$ be a continuous $*$-linear mapping. Then every $\sigma$-derivation $d$ is automatically continuous.\end{theorem}

\begin{proof} By Proposition 3.4 we may assume that $\sigma$ is a
continuous $*$-homomorphism. Theorem 2.1 implies that $d$ is an
ordinary derivation from ${\cala}$ into $\widetilde{B({\calh})}$.
By the main theorem of \cite{RIN} we conclude that $d$ is
continuous. \end{proof}

Using the same argument as in the proof of Proposition 3.6 we can
establish the following theorem:

\begin{theorem} Let $\sigma$ be a continuous linear mapping. Then every $*$-$\sigma$-derivation is automatically continuous.
\end{theorem}

\section{$(\sigma,\tau)$-derivations on $C^*$-algebras}

In \cite{M-M} the authors considered the notion of
$(\sigma,\tau)$-derivation. Assume that ${\cala}$ is a
$*$-subalgebra of a $*$-algebra ${\calb}$ and
$\sigma,\tau:{\cala}\to {\calb}$ are $*$-linear mappings. A linear
mapping $d:{\cala}\to {\calb}$ is called a
$*$-$(\sigma,\tau)$-derivation if $d$ preserves $*$ and
$d(ab)=d(a)\sigma(b)+\tau(a)d(b)$ for all $a,b\in {\cala}$.
Obviously, a $*$-$\sigma$-derivation is a
$*$-$(\sigma,\sigma)$-derivation.

\begin{theorem} Let ${\cala}$ be a $*$-subalgebra of a $*$-algebra ${\calb}$ and $\sigma,\tau:{\cala}\to {\calb}$ are $*$-linear mappings. Then every $*$-$(\sigma,\tau)$-derivation $d:{\cala}\to {\calb}$ is a $*$-$(\frac{\sigma+\tau}{2},\frac{\sigma+\tau}{2})$-derivation.\end{theorem}

\begin{proof} First we show that each $*$-$(\sigma,\tau)$-derivation is
a $*$-$(\tau,\sigma)$-derivation. We have
$d(ab)=d(b^*a^*)^*=(d(b)^*\sigma(a)^*+\tau(b)^*d(a)^*)^*=d(a)\tau(b)+\sigma(a)d(b)$.

Now we conclude that
$d(ab)=\frac{1}{2}d(ab)+\frac{1}{2}d(ab)=\frac{1}{2}(d(a)\sigma(b)+\tau(a)d(b))+\frac{1}{2}(d(a)\tau(b)+\sigma(a)d(b))=d(a)\frac{\tau+\sigma}{2}(b)+\frac{\tau+\sigma}{2}(a)d(b)$.
\end{proof}

The previous theorem enables us to focus on $\sigma$-derivations
while we deal with $*$-algebras. In particular, by using Theorem
4.1 we obtain the following generalizations of Theorem 3.2 and
Theorem 3.7.

\begin{theorem} Let $\sigma$ and $\tau$ be continuous $*$-linear mappings from a $C^*$-algebra ${\cala}$ acting on a Hilbert space ${\calh}$ into $B({\calh})$ and let $d:{\cala}\to B({\calh})$ be a continuous $*$-$(\sigma,\tau)$-derivation. Then there is a continuous linear mapping $\Sigma:{\cala}\to B({\calh})$ such that $d$ is a $*$-$\Sigma$-derivation.\end{theorem}

\begin{theorem} If $\sigma$ and $\tau$ are continuous $*$-linear mappings from a $C^*$-algebra ${\cala}$ acting on a Hilbert space ${\calh}$ into $B({\calh})$, then every $*$-$(\sigma,\tau)$-derivation $d:{\cala}\to B({\calh})$ is automatically continuous.\end{theorem}

\bibliographystyle{amsplain}

\end{document}